\documentclass[10pt]{article}

\usepackage{latexsym}
\usepackage{epsfig}
\usepackage{amssymb,amsfonts,amscd}
\usepackage{amsmath}

\newcommand{\bC}{{\mathbb C}}
\newcommand{\dt}{{\Delta t}}
\newcommand{\dz}{{\Delta z}}
\newcommand{\ubar}{{\bar u}}
\newcommand{\vbar}{{\bar v}}
\newcommand{\K}{\kappa}
\newcommand{\D}{\Delta}
\newcommand{\Dbar}{{\bar\D}}

\newcommand{\Jinv}{\|J^{-1}\|}
\newcommand{\HOT}{\hbox{\em H.O.T.}}
\newcommand{\tol}{\tau}
\newcommand{\Pbar}{{P'}}
\newcommand\Cp{C$'$}

\newcommand{\bP}{{\mathbb P}}

\newcommand\qed{{\hspace*{\fill}$\Box$\vskip12pt plus 1pt}}

\newcommand\sE{{\mathcal E}}

\newcommand\sI{{\mathcal I}}

\newcommand\sO{{\mathcal O}}

\def\cond{{\mathop{\rm cond}\nolimits}}

\newtheorem{theorem}{Theorem}[section]

\newtheorem{rem}[theorem]{Remark}

\newtheorem{Example}[theorem]{Example}

\begin{document}

\title{Multiprecision path tracking}

\author{D. J. Bates\thanks{Department of Mathematics,
University of Notre Dame, Notre Dame, IN 46556 (dbates1@nd.edu,
http://www.nd.edu/$\sim$dbates1) This author was supported by the
Duncan Chair of the University of Notre Dame, the University of
Notre Dame, NSF grant DMS-0410047 and the Arthur J. Schmitt
Foundation}, A. J. Sommese\thanks{Department of Mathematics,
University of Notre Dame, Notre Dame, IN 46556 (sommese@nd.edu,
http://www.nd.edu/$\sim$sommese). This author was supported by the
Duncan Chair of the University of Notre Dame, the University of
Notre Dame, and NSF grant DMS-0410047}, and C. W.
Wampler\thanks{General Motors Research and Development, Mail Code
480-106-359, 30500 Mound Road, Warren, MI 48090
(Charles.W.Wampler@gm.com, http://www.nd.edu/$\sim$cwample1)  This
author was supported by NSF grant DMS-0410047}}

\date{May 3, 2006}
\maketitle

\begin{abstract}
A path tracking algorithm that adaptively adjusts precision is
presented. By adjusting the level of precision in accordance with
the numerical conditioning of the path, the algorithm achieves high
reliability with less computational cost than would be incurred by
raising precision across the board. We develop simple rules for
adjusting precision and show how to integrate these into an
algorithm that also adaptively adjusts the step size. The behavior
of the method is illustrated on several examples arising as
homotopies for solving systems of polynomial equations.



\end{abstract}

Homotopy continuation, numerical algebraic geometry, polynomial
systems

\section{Introduction}
Path tracking is the task of tracing out a 1 real dimensional
solution curve described implicitly by a system of equations,
typically $n$ equations in $n+1$ variables, given an initial point
on, or close to, the path.  This can arise in many ways, but our
motivation is the solution of systems of polynomials via homotopy
continuation (see \cite{AG90,Li97,Li03,Mor87,MS87,SoWa}). In this
method, to find the isolated solutions of the system $f(z)=0$ for
given polynomials $f:\mathbb C^n \to \mathbb C^n$, one constructs a
homotopy, $H(z,t)$, $H: \mathbb C^n \times \mathbb C \to \mathbb
C^n$ such that $H(z,0)=f(z)$ is the target system to be solved while
$H(z,1)$ is a starting system whose isolated solutions are known.
There is a well-developed theory on how to construct such homotopies
to guarantee, with probability one, that every isolated solution of
$f(z)=0$ is the endpoint in the limit as $t\rightarrow0$ of at least
one smooth path $z_i(t)$, where $H(z_i(t),t)=0$ on $t\in(0,1]$ and
where $z_i(1)$, $i=1,2,\ldots$, are the known isolated solutions of
$H(z,1)=0$.  Similar constructions arise in other contexts, where
the existence of a path leading to the desired solutions may or may
not be guaranteed.  Even when there is no guarantee, experience
shows that in some application domains continuation techniques yield
solutions more reliably than Newton's method, especially when good
initial guesses are not available.  While in our applications the
path is just a means of arriving at the endpoint, in other
applications one may desire to accurately trace out the path itself,
such as when plotting the response of a mathematical model as one of
its parameters is varied.

The most common path tracking algorithms are predictor-corrector
methods: from an approximate solution point on the path, a predictor
gives a new approximate point a given step size along the path, then
a corrector brings this new point closer to the path.  For example,
one may use an Euler predictor, which steps ahead along the tangent
to the path, or a higher order predictor that uses several recent
points and the derivatives of the homotopy function at them to
extrapolate to the predicted point.  Typically, the prediction is
then used as the initial point for correction by Newton's method.
Since the solution set is one-dimensional, an extra constraint is
introduced to isolate the target of the correction.  For general
homotopies, a useful constraint is to find where the solution path
intersects the hyperplane normal to the last computed tangent
direction. In the more restrictive setting of polynomial systems,
the homotopy can be designed such that the paths advance
monotonically with $t$, that is, there are no turning points, in
which case it is acceptable (and simpler) to perform corrections by
holding $t$ fixed.  The adaptive precision algorithm we discuss here
is compatible with any of these prediction and correction schemes.

For good results, the predictor step size must be chosen
appropriately. Too large a step size may result in a prediction
outside the zone of convergence of the corrector, while too small a
step size means progress is slow and costly. Consequently, it has
long been recognized that adaptive control of the step size is
crucial for obtaining good reliability without undue computational
cost.

While step size control is well established, less attention has been
paid to efficient handling of precision.  With wider availability of
software packages for higher precision arithmetic, along with faster
computers to execute the software, it becomes interesting to
consider how adjustable precision might be deployed to improve the
performance of path tracking algorithms.  The issue at stake is
analogous to step size control: without enough precision, path
tracking will fail, but the use of excessive precision is
inefficient.  To address this tradeoff, this paper proposes an
algorithm that dynamically adjusts the number of digits used in
computations according to the evolution of the numerical
conditioning of the homotopy function.

In our primary application of interest, the solution of polynomial
systems, there are several factors driving the need for higher
precision.  It is well known that high degree polynomials often lead
to ill-conditioned problems.  When treating polynomial systems in
several variables, the total degree of the system, being the product
of the degrees of the individual equations, quickly becomes large
even for low degree polynomials, which can also lead to
ill-conditioning.  Thus, one driving force is the desire to solve
larger systems of higher total degree.  A second motivation is that
our systems often have some, or possibly many, singular solutions,
and thus, the solution paths leading to these solutions are
necessarily ill-conditioned near the end. While endgame methods
exist for enhancing the accuracy with which such endpoints can be
estimated, for singularities of high enough multiplicity more
precision is required.  Finally, although the homotopy constructions
guarantee, with probability one, that no path passes exactly through
a singularity before reaching its endpoint, there is always a chance
that a near singular condition can be encountered.  To obtain the
highest reliability possible, we need to detect this and allocate
sufficient digits to successfully track past such obstructions.

The paper is organized as follows. In section 2, we review the
behavior of Newton's method in floating point, revealing how its
accuracy and convergence properties depend on precision. In section
3, we discuss path tracking with adaptive step size control and
identify how it fails when precision is insufficient. This leads, in
section 4, to a novel technique for path tracking using adaptive
precision. This new adaptive precision path tracking algorithm has
been implemented in a software package, Bertini, currently under
development by the authors. Several examples are presented in
section 5 to illustrate the usefulness of adaptive precision.
Finally, in section 6, a few related ideas that would make
interesting studies are discussed.

\section{Background: Newton's method}\label{Sec:Newton}

The core numerical process in the path tracker is the corrector,
which in our case is Newton's method.  A good predictor speeds up
the path tracker by allowing a large step while still supplying an
initial guess within the convergence region of the corrector.
However, it is the loss of convergence that causes path tracking to
fail. In exact arithmetic, as long as the path remains nonsingular,
there must be a region surrounding the path within which Newton's
method converges quadratically. With a small enough step $\dt$ in
$t$, we can be assured of advancing along the path, although
possibly very slowly.  This holds even if we use only a zero-th
order predictor, i.e., if the point from the last value $t_k$ is
used to initialize the corrector for the new value
$t_{k+1}=t_k+\dt$. In contrast, in inexact floating point
arithmetic, the convergence region can disappear, thus halting the
path tracker.  Short of this, an unacceptably slow linear rate of
convergence might dominate, causing the step size to plummet.  It
can also happen that the corrector converges but to an answer that
is outside the desired tolerance.

Due to these considerations, an analysis of Newton's method in
floating point is of interest and will help us derive rules for
setting the precision used in our path tracker.  The following
analysis resembles that of \cite{Tisseur}. Let
$F(z):\bC^n\rightarrow\bC^n$ be continuously differentiable, and
denote its Jacobian matrix of partial derivatives as $J(z)$.  To
solve $F(z)=0$ by Newton's method given an initial guess $z_0$, one
iteratively updates $z_i$, $i=1,2,\dots$, as
\begin{equation}\label{Eq:Newton}
    \begin{array}{l}
       \hbox{Solve } J(z_i)\dz_i = -F(z_i) \hbox{ for } \dz_i, \\
       z_{i+1} = z_i + \dz_i.
     \end{array}
\end{equation}

Suppose that we work in floating point with unit roundoff $u$.  In
other words, if we compute with a mantissa of $b$ bits in binary or
with $P$ digits in decimal, $u=2^{-b}=10^{-P}$.  We may consider
evaluating the residuals $F(z_i)$ in higher precision, say $\ubar\le
u$. Let $\hat F(z)$ be the floating point output of the procedure
that evaluates $F(z)$.  We assume that there exists a function
$\psi$ depending on $z$, $u$, and $\ubar$ such that the error
$e(z)=\hat F(z)-F(z)$ obeys
\begin{equation}\label{Eq:errorF}
    \| e(z) \| \le u\|F(z)\| + \psi(z,u,\ubar).
\end{equation}
By definition, at a solution point $z_*$, we have $F(z_*)=0$, so it
is clear that the function $\psi$ drives the final error.  To
determine $\psi$, one must examine the function $F$ and the program
that implements it.  We will give a rough rule of thumb later for
the systems we treat.

In solving Eq.~\ref{Eq:Newton} for the correction $\dz_i$, there is
error in evaluating $J(z_i)$ and in solving the linear system.  Both
errors can be absorbed into an error matrix $E_i$ such that the
computed correction is
\begin{equation}\label{Eq:errorEdef}
    \dz_i = (J(z_i)+E_i)^{-1}(F(z_i)+e(z_i)).
\end{equation}
We assume this error is bounded by
\begin{equation}\label{Eq:errorEbound}
    \|E_i\| \le \sE \left( u \|J(z_i)\|+\phi(z_i,u)\right),
\end{equation}
for some constant $\sE>1$ and positive function $\phi$.  We expect
the first term because of roundoff of the Jacobian, whereas $\phi$
accounts for errors in evaluating $J$ that do not vanish even when
$J$ does.  The constant $\sE$ accounts for the subsequent growth in
the error during the linear solve.

For simplicity of notation, let $v=z_i$ be the current guess,
$\vbar=z_{i+1}$ the new guess after a single iteration, and let
$v_*$ be the solution point near $v$.  Also, let's use the shorthand
notations $F=F(v)$, $J=J(v)$, $J_*=J(v_*)$, $\D=\|v-v_*\|$ and
$\Dbar=\|\vbar-v_*\|$. In the next paragraph, we will establish a
bound on $\Dbar$ in terms of $\D$.  Whenever $\Dbar<\D$, the Newton
step successfully reduces the error in the estimate of the root
$v_*$.

Since $F(v_*)=0$, the Taylor series of $F(z)$ at $v_*$ gives
\[
  F(z) = J_*\cdot(z-v_*) + \HOT
\]
where the higher order terms, $\HOT$, are quadratic or higher in
$z-v_*$.  Similarly,
\[
  J(z) = J_* + \HOT
\]
where the higher order terms are linear in $z-v_*$. Consequently, in
a ball $B=\{z\>:\>\|z-v_*\|\le R\}$ centered on $v_*$ with $v\in B$,
there exist positive constants $\alpha$ and $\beta$ such that
\begin{align}
  &\|F(z)\| \le \|J_*\|\|z-v_*\| + \alpha\|z-v_*\|^2,\label{Eq:alpha1}\\
  &\|F(z)-J_*(z-v_*)\| \le \alpha\|z-v_*\|^2\label{Eq:alpha2}\\
  &\|J_*\|\le\|J\|+\beta\|z-v_*\|, \qquad \|J-J_*\|\le \beta\|z-v_*\|.\label{Eq:beta1}
\end{align}
In Newton's method, we solve
\begin{equation}\label{Eq:NewtonStepWithError}
  (J+E)d=-(F(v)+e)
\end{equation}
for $d$ and take the step
\begin{equation}\label{Eq:vbarUpdate}
  \vbar=v+d+\varepsilon,
\end{equation}
where $\varepsilon$ is the error in forming the sum.  The standard
model of round-off error in floating point addition \cite{Wilk}
gives
\begin{equation}\label{Eq:AdditionRoundoff}
 \|\varepsilon\| \le u(\|v\|+\|d\|) \le u(\D + \|v_*\| + \|d\|),
\end{equation}
so subtracting $v_*$ from both sides of Eq.~\ref{Eq:vbarUpdate}, we
have
\begin{equation}\label{Eq:Dbar1}
  \Dbar \le \|v-v_*+d\| + u(\D + \|v_*\| + \|d\|).
\end{equation}
If $J$ is nonsingular and $\Jinv\|E\|<1$, then $(J+E)$ is
nonsingular, the Newton step is well defined, and
\begin{equation}\label{Eq:KJinv}
  \|(J+E)^{-1}\| \le K\Jinv, \qquad K=\frac{1}{1-\Jinv\|E\|}.
\end{equation}
Accordingly, from Eq.~\ref{Eq:NewtonStepWithError} we have
\begin{equation}\label{Eq:dnorm}
    \|d\|\le K \Jinv (\|F\| + \|e\|).
\end{equation}
Also, after adding $(J+E)(v-v_*)$ to both sides of
Eq.~\ref{Eq:NewtonStepWithError} and simplifying using
Eqs.~\ref{Eq:alpha1}--\ref{Eq:beta1}, we have
\begin{equation}\label{Eq:dvdnorm}
    \|v-v_*+d\| \le K\Jinv( \|E\|\D + (\alpha+\beta)\D^2+\|e\| ).
\end{equation}
Substituting from Eqs.~\ref{Eq:dnorm},\ref{Eq:dvdnorm} into
Eq.~\ref{Eq:Dbar1} and using
Eqs.~\ref{Eq:errorF},\ref{Eq:errorEbound}, we obtain the bound
\begin{align}\label{Eq:Dbar}
    \Dbar \le & K\Jinv(1+u)^2(\alpha+\beta)\D^2 + \nonumber\\
      &\left(K\Jinv[(2+\sE+u)\left(u\|J\|+\phi\right)]+u \right)\D + K\Jinv(1+ u)\psi + u\|v_*\|.
\end{align}
This relation holds as long as $\Jinv\|E\|<1$, so that the linear
solve in the Newton step is well-defined, and $v$ is in the ball
$B$, so that Eqs.~(\ref{Eq:alpha1}--\ref{Eq:beta1}) hold.

If $\Dbar<\D$, the Newton step reduces the error in the
approximation of the root. In exact arithmetic, we have
$u=\phi=\psi=0$ and $K=1$, so $\Dbar\le \Jinv(\alpha+\beta)\D^2$.
The error contracts if the initial guess is accurate enough so that
$\D<1/(\Jinv(\alpha+\beta))$.  If we also have
$\D<1/(\|J_*^{-1}\|(\alpha+\beta))$, it is clear that all subsequent
iterates are nonsingular and contractive, from which one has the
well-known result that Newton's method converges quadratically to a
nonsingular solution for a sufficiently accurate initial guess.  One
sees that the more singular the Jacobian is at the root, the slower
the convergence and the smaller the convergence zone.

In floating point arithmetic, we cannot expect the error to converge
to zero.  From Eq.~\ref{Eq:Dbar}, one may expect the error to
contract until
\begin{equation}\label{Eq:FinalError}
    \D \approx (1+u)K\Jinv\psi + u\|v_*\| \approx K\|J_*^{-1}\|\psi +
    u\|v_*\|.
\end{equation}
The second term is the error inherent in representing $v_*$ in
floating point.  The first term depends on the accuracy, $\psi$,
with which the function is evaluated.  This can be reduced by using
higher precision, $\ubar$, in the function evaluation per
Eq.~\ref{Eq:errorF}. The precision of the Jacobian and the linear
solve do not affect the final error.

On the other hand, the precision of the Jacobian and the linear
solve do affect convergence.  Without enough precision, $\Jinv\|E\|$
may approach or surpass 1, which means that the linear solve may
fail due to singularity or may yield such an inaccurate step that
the error diverges.  Notice that $\Jinv\|E\| = u\Jinv\|J\| +
\Jinv\phi = u\K + \Jinv\phi$, where $\K=\cond(J)$.  The first term,
$u\K$, reflects the well-known result that in solving a linear
system, floating point roundoff is magnified by the condition number
of the matrix.

\section{Step length control and failure}

To produce an improved path tracking algorithm, it is useful to
first examine a standard predictor/corrector algorithm to see why
adaptive step length control generally succeeds when conditioning is
mild and why it may fail when conditioning is adverse.

A simple and effective approach for step length control is to adjust
the step length up or down according to the success or failure of a
complete prediction/correction cycle.  Suppose the homotopy function
$H(z,t)=0$ defines a one-dimensional nonsingular path $z(t)$. We are
given a start point approximately on the path, $z_0\approx z(t_0)$,
an ending value of $t$, and a tolerance to which points on the path
are to be found. Then, in brief, a predictor/corrector path tracker
with adaptive step length control may be constructed as follows.
\begin{description}
    \item[Initialize] Select: an initial step size, $s$; the
    number of corrector iterations allowed per step, $N\ge1$; the step
    adjustment factor, $a\in(0,1)$; the step expansion integer, $M\ge1$; and a minimum
    step size $s_{\rm min}$.
    \item[Predict] Estimate a new point near the path whose distance from the
    current point is the step size.
    \item[Correct] Iteratively improve the new path point,
    constraining its distance from the prior path point.  Allow at
    most $N$ iterations to reach the specified tolerance.
    \item[On success] If the tolerance is achieved:
    \begin{itemize}
        \item Update the current path point to be the newly found path point.
        \item If we have reached the final value of $t$, exit with
        \emph{success}.
        \item If there have been $M$ successes in a row, expand the step
        size by $s=s/a$.
    \end{itemize}
    \item[On failure] If the tolerance is not achieved:
    \begin{itemize}
        \item Cut the step length by $s=as$.
        \item If $s<s_{\rm min}$, exit with \emph{failure}.
    \end{itemize}
    \item[Loop] Go back to {\bf Predict}.
\end{description}

The key is to allow only a small number of iterations in the
corrector, typically only $N=2$ or $N=3$. This forces the prediction
to stay within a good convergence region surrounding the path. If a
large number of iterations is allowed, a bad prediction might
ultimately converge, but it may wander first and become attracted to
a completely different path in the homotopy. Keeping $N$ small, the
step size adaptation slows down to negotiate sharp turns in the path
and accelerates whenever the path is relatively straight. Properly
implemented, this results in a robust and efficient path tracking
algorithm.

We can be a bit more precise.  Let us do so by specifically
considering an Euler predictor with a Newton corrector.  Both of
these derive from the linearized local model of the path. The Taylor
series at $(z_1,t_1)$ is
\begin{equation}\label{Eq:TaylorSeries}
  H(z_1+\dz,t_1+\dt) = H(z_1,t_1) +
     \frac{\partial H}{\partial z}(z_1,t_1)\dz +
     \frac{\partial H}{\partial t}(z_1,t_1)\dt +
     \HOT
\end{equation}
where the higher order terms, $\HOT$, are quadratic or higher in
$(\dz,\dt)$.  Ignoring the higher order terms and setting
$H(z_1+\dz,t_1+\dt)=0$, one has the basic Euler predictor and Newton
corrector relations.  These are a system of $n$ equations in $n+1$
unknowns; as long as the combined matrix $[\partial H/\partial
z\>\>\partial H/\partial t]$ is rank $n$, there is a well-defined
tangent direction and tracking may proceed.  The predictor adds a
constraint on the length of the step along the tangent, whereas
corrector steps are constrained to move transverse to the tangent.
The extra constraints are particularly simple in the case where
$\partial H/\partial z$ is rank $n$, for then the path progresses
monotonically in $t$, and the step can be controlled via the advance
of $t$.  Accordingly, one has a linear system to be solved for
$\dz$:
\begin{equation}\label{Eq:BasicStep}
    \left[\frac{\partial H}{\partial z}(z_1,t_1)\right]\dz =
      -\left( H(z_1,t_1) +
        \frac{\partial H}{\partial t}(z_1,t_1)\dt \right).
\end{equation}
For prediction, we set $\dt=s$, the current step size, and for
correction, we set $\dt=0$.

Since the neglected terms are quadratic, the prediction error is
order $\sO(s^2)$.  Thus, in the case of a failed step, cutting the
step size from $s$ to $as$ reduces the prediction error by a factor
of $a^2$.  In this way, cuts in the step size quickly reduce the
prediction error until it is within the convergence region of the
corrector.  With a $k$th order predictor, the prediction error
scales as $a^{k+1}$, potentially allowing larger step sizes.  In any
case, the adaptive approach quickly settles to a step size $s$ just
small enough so that the corrector converges, while the next larger
step of $s/a$ fails.  With $a=1/2$ and $M=5$, the step size adapts
to within a factor of 2 of its optimum, with an approximate overhead
of 20\% spent checking if a larger step size is feasible.

Failure of path tracking with an adaptive step size can be
understood from the discussion of Newton's method in
\S~\ref{Sec:Newton}.  For small enough initial error and
infinite-precision arithmetic, the Newton corrector gives quadratic
convergence to a nonsingular root. Near a singularity,
$\|J_*^{-1}\|$ is large, which can lead to a small quadratic
convergence zone and a slower rate of quadratic convergence. Inexact
arithmetic can further shrink the convergence zone, degrade the
convergence rate from quadratic to linear, and introduce error into
the final answer. From these considerations, we see that there are
two ways for the adaptive step size path tracker to halt prematurely
near a singularity.
\begin{enumerate}
    \item The predictor is limited to a tiny
    step size to keep the initial guess within the convergence zone
    of the corrector.  If this is too small, we may exceed the
    allotted computation time for the path.
    \item The path may approach a point where the final error of
    the corrector is as large as the requested path tracking
    tolerance.
\end{enumerate}

The first mode of failure can occur even with infinite precision,
but degradation of the convergence properties with too low a
precision increases the occurrence of this failure.  The second mode
of failure is entirely a consequence of lack of precision. By
allocating enough precision, we can eliminate the second mode of
failure and reduce the occurrence of the first mode.  It is
important to note that in some applications there is flexibility in
the definition of the homotopy, which can be used to enlarge
convergence zones and thereby speed up path tracking.  For example,
re-scaling of the equations and variables can sometimes help.
However, such maneuvers are beyond the scope of this paper, which
concentrates only on tracking the path of a given homotopy.

\section{Adaptive Precision}

The use of high precision can largely eliminate both types of path
tracking failure identified above.  However, high precision
arithmetic is expensive, so it must be employed judiciously. One
might be tempted to rachet precision up or down in response to step
failures as in the adaptive step size algorithm.  This presents the
difficulty that there is just one stimulus, step failure, and two
possible responses, cut the step size or increase precision.  In the
following paragraphs, we outline two possible algorithms for
adapting both step size and precision.

\subsection{Adapting precision via path re-runs}

The simplest approach to adapting precision, shown in Figure 1, is
to run the entire path in a fixed precision with adaptive re-runs.
That is, if the path tracking fails, one re-runs it in successively
higher precision until the whole path is tracked successfully or
until limits in computing resources force termination.  The
advantage of this approach is that adaptation is completely external
to the core path tracking routine. Thus, this strategy can be
applied to any path tracker that enables requests for higher
precision.  For example, in the polynomial domain, the package PHC
\cite{V99} offers multiple precision, although the precision must be
set when calling the program.

The adaptation algorithm of Figure~1 has two main disadvantages.
First, when too low a precision is specified, the tracker may waste
a lot of computation near the point of failure before giving up and
initiating a re-run in higher precision. Second, the whole path is
computed in high precision when it may be needed only in a small
section of the path, often near the end in the approach to a
singular solution. A slightly more sophisticated treatment can avoid
re-computing the segment of the path leading up to the failure point
by requesting the tracker to return its last successful path point.
The re-run in higher precision can then be initiated from that point
on.

\begin{figure}\label{Fig:fixed_flow}
\begin{picture}(14.000000,11.000000)(0.000000,-11.000000)
\put(3.0000,-1.0000){\oval(6.0000,2.0000)}
\put(0.0000,-2.0000){\makebox(6.0000,2.0000)[c]{\shortstack[c]{
Homotopy, start points,\\
initial constant settings }}}
\put(3.0000,-2.0000){\vector(0,-1){1.0000}}
\put(1.0000,-5.0000){\framebox(4.0000,2.0000)[c]{\shortstack[c]{
Path following\\
with adaptive\\
stepsize }}} \put(3.0000,-5.0000){\vector(0,-1){1.0000}}
\put(1.6250,-7.3750){\line(1,1){1.3750}}
\put(1.6250,-7.3750){\line(1,-1){1.3750}}
\put(4.3750,-7.3750){\line(-1,-1){1.3750}}
\put(4.3750,-7.3750){\line(-1,1){1.3750}}
\put(1.6250,-8.7500){\makebox(2.7500,2.7500)[c]{\shortstack[c]{
Path\\
success? }}} \put(4.3750,-6.9625){\makebox(0,0)[lt]{No}}
\put(3.4125,-8.7500){\makebox(0,0)[lb]{Yes}}
\put(3.0000,-8.7500){\vector(0,-1){1.0000}}
\put(1.5000,-11.0000){\framebox(3.0000,1.2500)[c]{\shortstack[c]{
Return\\
endpoint }}} \put(4.3750,-7.3750){\line(1,0){4.0000}}
\put(8.3750,-7.3750){\vector(0,1){1.0000}}
\put(6.7500,-4.7500){\line(1,1){1.6250}}
\put(6.7500,-4.7500){\line(1,-1){1.6250}}
\put(10.0000,-4.7500){\line(-1,-1){1.6250}}
\put(10.0000,-4.7500){\line(-1,1){1.6250}}
\put(6.7500,-6.3750){\makebox(3.2500,3.2500)[c]{\shortstack[c]{
Reruns\\
remaining? }}} \put(8.8625,-3.1250){\makebox(0,0)[lt]{Y}}
\put(10.0000,-4.2625){\makebox(0,0)[lt]{N}}
\put(10.0000,-4.7500){\vector(1,0){1.0000}}
\put(11.0000,-5.3750){\framebox(3.0000,1.2500)[c]{\shortstack[c]{
Failed\\
Path }}} \put(8.3750,-3.1250){\line(0,1){0.5000}}
\put(8.3750,-2.6250){\vector(0,1){1.0000}}
\put(6.8750,-1.6250){\framebox(3.0000,1.2500)[c]{\shortstack[c]{
Increase\\
precision }}} \put(6.8750,-1.0000){\vector(-1,0){0.8750}}
\end{picture}
\caption{Adapting precision via path re-runs}
\end{figure}
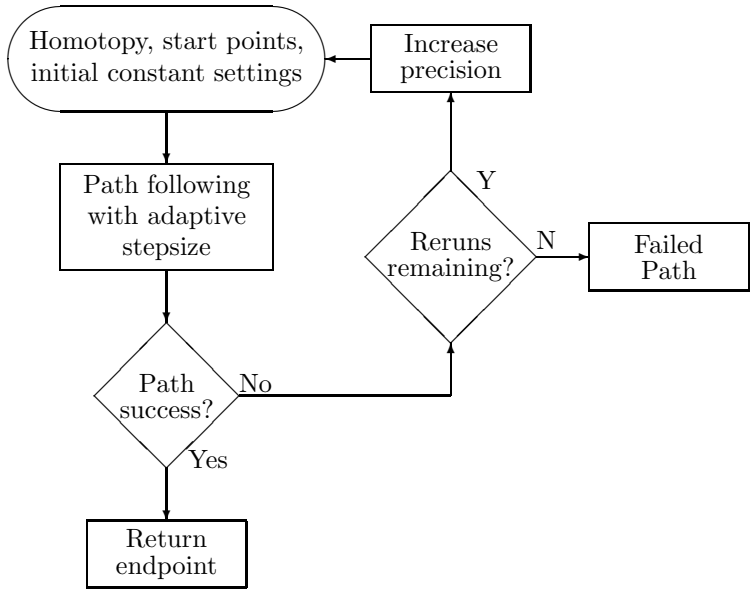

\subsection{Stepwise adaptive precision}\label{Sec:Step_adaptive}

Instead of waiting for the adaptive step size method to fail before
initiating higher precision, we propose to continuously monitor the
conditioning of the homotopy to judge the level of precision needed
at each step.  In this way, the computational burden of higher
precision is incurred only as needed, adjusting up and down as the
tracker proceeds, while obtaining superior reliability.

To decide how much precision is needed, we turn to the analysis of
Newton's method from \S~\ref{Sec:Newton}.  We wish to ensure that
the achievable accuracy is within the specified tolerance and that
convergence is fast enough.

In what follows, we need to evaluate $\|J\|$ and $\Jinv$.  These do
not need to be very accurate, as we will always include safety
margins in the formulas that use them.  $\|J\|$ is readily available
in the max norm, where we use the maximum magnitude of any of its
entries.  $\Jinv$ is more difficult, as we do not wish to compute
the full inverse of the matrix.  This issue has been widely studied
in terms of estimating the condition number $\K=\|J\|\Jinv$.  A
relatively inexpensive method, suggested in \cite{Wat} and
elsewhere, is to choose a unit vector $b$ at random and solve $Jy=b$
for $y$. Then, we use the estimate $\Jinv\approx\|y\|$. Although
this underestimates $\Jinv$, tests of matrices up to size $10\times
10$ show the approximation to be
reliably within a factor of 10 of the true value, which is easily
absorbed into our safety margins.

One requirement is that $\Jinv\|E\|$ should be small enough to
ensure that the error-perturbed Jacobian is nonsingular.  Minimally,
we require $\Jinv\|E\|<1$, but by requiring it to be a bit smaller,
say $\Jinv\|E\|<10^{-\sigma_1}$ for some $\sigma\ge1$, we force
$K\approx 1$.  This removes the growth of $K$ as one possible source
of failure.  Suppose that the error function $\phi$ in
Eq.~\ref{Eq:errorEbound} is of the form $\phi=\Phi u$.  Then, our
first rule is to require
\begin{equation}\label{Eq:setK}
    \Jinv\sE\left(\|J\| + \Phi\right)u < 10^{-\sigma_1}.
\end{equation}
Using $P$ decimal digits of arithmetic results in precision
$u=10^{-P}$, so we may restate this rule as
\renewcommand{\theequation}{{\bf A}}
\begin{equation}\label{Eq:setKlog}
    P > \sigma_1+ \log_{10}[\Jinv\sE\left(\|J\| + \Phi\right)].
\end{equation}
\renewcommand{\theequation}{\arabic{equation}}

A second requirement is that the corrector must converge within $N$
iterations, where we keep $N$ small as in the usual adaptive step
size algorithm, typically 2 or 3. Let us say that the tolerance for
convergence is $\D=\|v-v_*\|<10^{-\tol}$.  Recall that in each step
of Newton's method, we compute $d$ and take the step $\vbar=v+d$.
The best estimate available of the accuracy is $\D\approx\|d\|$, so
we declare success when $\|d\|<10^{-\tol}$.  Suppose that after
$i<N$ iterations this is not yet satisfied.  We still have $N-i$
iterations to meet the tolerance, and we would like to be sure that
a lack of precision does not prevent success. Pessimistically, we
assume that the linear factor in $\D$ in Eq.~\ref{Eq:Dbar} dominates
the quadratic one and that the rate of convergence does not improve
with subsequent iterations. We force $K\approx 1$, and we have $u\ll
1$. Including the same safety margin as before, $10^{\sigma_1}$, the
requirement becomes
\begin{equation}\label{Eq:setConvergence}
    \left(10^{\sigma_1}\Jinv(2+\sE)(u\|J\|+\phi)+u \right)^{N-i}\|d\|<10^{-\tol}.
\end{equation}
As before, let's assume $\phi=\Phi u$.  Taking logarithms, the
number of decimal digits of precision must satisfy
\renewcommand{\theequation}{{\bf B}}
\begin{equation}\label{Eq:setConvergenceLog}
    P > \sigma_1 + \log_{10} \left(\Jinv(2+\sE)(\|J\|+\Phi)+1 \right) + (\tol+\log_{10}\|d\|)/(N-i).
\end{equation}
\renewcommand{\theequation}{\arabic{equation}}%
Since we only apply this formula when the tolerance is not yet
satisfied, we have $\|d\|>10^{-\tau}$, or equivalently,
$\tol+\log_{10}\|d\|
> 0$. This implies that between corrector iterations,
requirement~\ref{Eq:setConvergenceLog} is always more stringent than
Eq.~\ref{Eq:setKlog}. However, we still use Eq.~\ref{Eq:setKlog}
outside the corrector, because $\|d\|$ is not then available.

Our third requirement is that the precision must be high enough to
ensure that the final accuracy of the corrector is within the
tolerance at full convergence.  For this, Eq.~\ref{Eq:FinalError} is
binding, so including a safety margin of $10^{-\sigma_2}$ and using
the norm of the current approximate solution, $\|v\|$, as the best
available estimate of $\|v_*\|$, we require
\begin{equation}\label{Eq:setFinalError}
    \|J^{-1}\|\psi + u\|v\| < 10^{-\tol-\sigma_2}.
\end{equation}
Suppose the error in evaluating the homotopy function is given by
$\psi=\Psi \ubar$.  If the function is evaluated in the same
precision as the rest of the calculations, i.e., $\ubar=u$, we have
the requirement
\renewcommand{\theequation}{{\bf C}}
\begin{equation}\label{Eq:setFinalErrorLog1}
    P > \sigma_2 + \tol + \log_{10}(\|J^{-1}\|\Psi + \|v\|).
\end{equation}
If instead we evaluate the function to higher precision, say
$\ubar=10^{-\Pbar}<u=10^{-P}$, we have the dual criteria
\renewcommand{\theequation}{{\bf \Cp}}
\begin{equation}\label{Eq:setFinalErrorLog2}
    P > \sigma_2 + \tol + \log_{10}\|v\|,
      \qquad \Pbar > \sigma_2 + \tol + \log_{10}\|J^{-1}\| + \log_{10}\Psi.
\end{equation}
\renewcommand{\theequation}{\arabic{equation}}
The effect of adding the two errors is absorbed into the safety
factor $\sigma_2$.

Conditions A, B, and C (or \Cp) allow one to adjust the precision as
necessary without waiting for the adaptive step size to fail.  If
necessary, the precision can even be increased between corrector
iterations.  An algorithm using these criteria is described by the
flowchart in Figure 2.  In this flowchart, ``Failure'' in the
predictor or corrector steps means that the linear solve of
Eq.~\ref{Eq:BasicStep} has aborted early due to singularity.  Using
the magnitude of the largest entry in $J$ as $\|J\|$, Gaussian
elimination with row pivoting may declare such a failure when the
magnitude of the largest available pivot is smaller than
$u\sE\|J\|$, for then the answer is meaningless.  This is more
efficient than completing the linear solve and checking condition A
or B, as these are sure to fail.

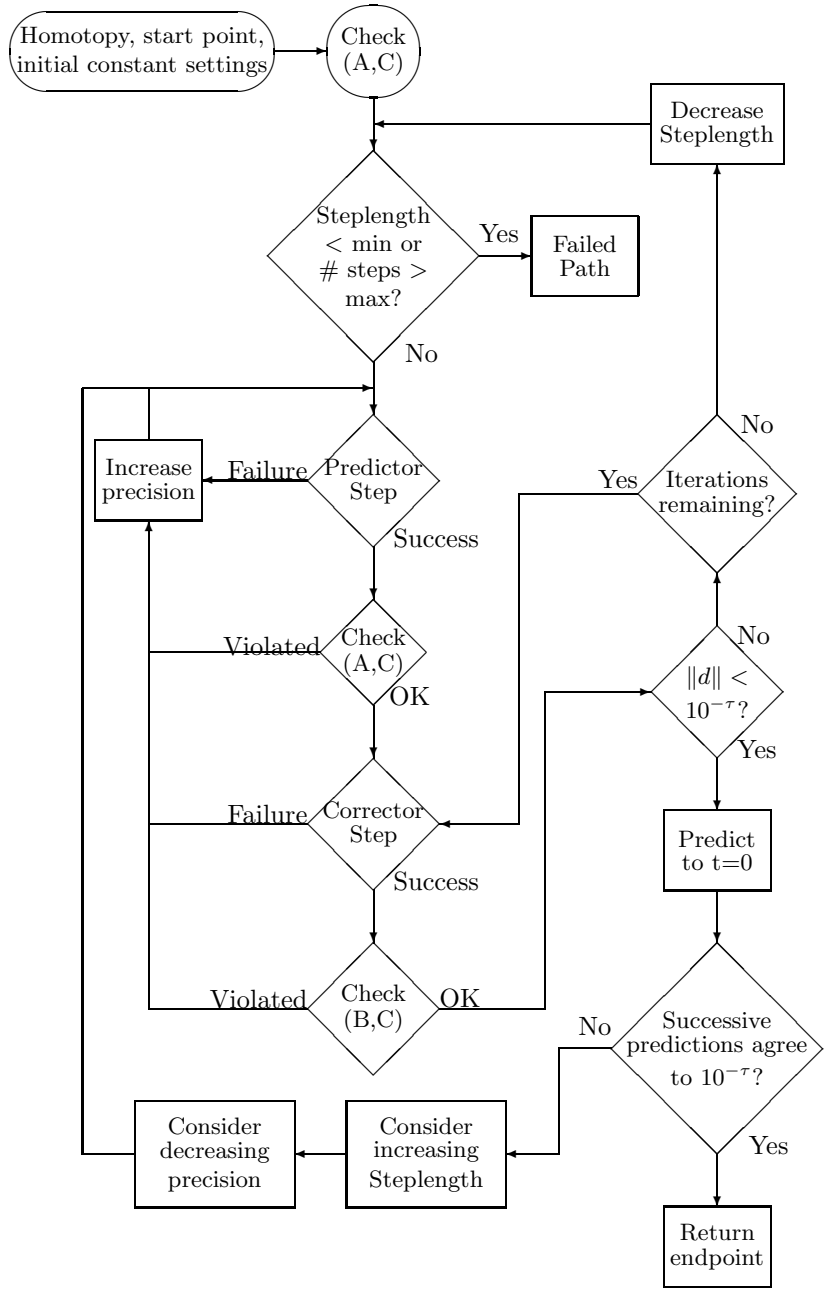
\begin{figure}\label{Fig:StepwiseFlowchart}
\begin{picture}(15.375000,24.250000)(0.000000,-24.125000)
\put(2.5000,-0.7500){\oval(5.0000,1.5000)}
\put(0.0000,-1.5000){\makebox(5.0000,1.5000)[c]{\shortstack[c]{
\small Homotopy, start point,\\
\small initial constant settings }}}
\put(5.0000,-0.7500){\vector(1,0){1.0000}}
\put(6.8750,-0.7500){\oval(1.7500,1.7500)}
\put(6.0000,-1.6250){\makebox(1.7500,1.7500)[c]{\shortstack[c]{
\small Check\\
\small (A,C) }}} \put(6.8750,-1.6250){\vector(0,-1){1.0000}}
\put(4.8750,-4.6250){\line(1,1){2.0000}}
\put(4.8750,-4.6250){\line(1,-1){2.0000}}
\put(8.8750,-4.6250){\line(-1,-1){2.0000}}
\put(8.8750,-4.6250){\line(-1,1){2.0000}}
\put(4.8750,-6.6250){\makebox(4.0000,4.0000)[c]{\shortstack[c]{
\small Steplength\\
\small $<$ min or\\
\small \# steps $>$ \\
\small max? }}} \put(8.8750,-4.0250){\makebox(0,0)[lt]{Yes}}
\put(7.4750,-6.6250){\makebox(0,0)[lb]{No}}
\put(8.8750,-4.6250){\vector(1,0){1.0000}}
\put(9.8750,-5.3750){\framebox(2.0000,1.5000)[c]{\shortstack[c]{
\small Failed\\
\small Path }}} \put(6.8750,-6.6250){\vector(0,-1){1.0000}}
\put(5.6250,-8.8750){\line(1,1){1.2500}}
\put(5.6250,-8.8750){\line(1,-1){1.2500}}
\put(8.1250,-8.8750){\line(-1,-1){1.2500}}
\put(8.1250,-8.8750){\line(-1,1){1.2500}}
\put(5.6250,-10.1250){\makebox(2.5000,2.5000)[c]{\shortstack[c]{
\small Predictor\\
\small Step }}} \put(5.6250,-8.5000){\makebox(0,0)[rt]{Failure}}
\put(7.2500,-10.1250){\makebox(0,0)[lb]{Success}}
\put(5.6250,-8.8750){\line(-1,0){1.0000}}
\put(4.6250,-8.8750){\vector(-1,0){1.0000}}
\put(1.6250,-9.6250){\framebox(2.0000,1.5000)[c]{\shortstack[c]{
\small Increase\\
\small precision }}} \put(2.6250,-8.1250){\line(0,1){1.0000}}
\put(2.6250,-7.1250){\vector(1,0){4.2500}}
\put(6.8750,-10.1250){\vector(0,-1){1.0000}}
\put(5.8750,-12.1250){\line(1,1){1.0000}}
\put(5.8750,-12.1250){\line(1,-1){1.0000}}
\put(7.8750,-12.1250){\line(-1,-1){1.0000}}
\put(7.8750,-12.1250){\line(-1,1){1.0000}}
\put(5.8750,-13.1250){\makebox(2.0000,2.0000)[c]{\shortstack[c]{
\small Check\\
\small (A,C) }}} \put(5.8750,-11.8250){\makebox(0,0)[rt]{Violated}}
\put(7.1750,-13.1250){\makebox(0,0)[lb]{OK}}
\put(5.8750,-12.1250){\line(-1,0){3.2500}}
\put(2.6250,-12.1250){\vector(0,1){2.5000}}
\put(6.8750,-13.1250){\vector(0,-1){1.0000}}
\put(5.6250,-15.3750){\line(1,1){1.2500}}
\put(5.6250,-15.3750){\line(1,-1){1.2500}}
\put(8.1250,-15.3750){\line(-1,-1){1.2500}}
\put(8.1250,-15.3750){\line(-1,1){1.2500}}
\put(5.6250,-16.6250){\makebox(2.5000,2.5000)[c]{\shortstack[c]{
\small Corrector\\
\small Step }}} \put(5.6250,-15.0000){\makebox(0,0)[rt]{Failure}}
\put(7.2500,-16.6250){\makebox(0,0)[lb]{Success}}
\put(5.6250,-15.3750){\line(-1,0){3.0000}}
\put(2.6250,-15.3750){\line(0,1){5.0000}}
\put(6.8750,-16.6250){\vector(0,-1){1.0000}}
\put(5.6250,-18.8750){\line(1,1){1.2500}}
\put(5.6250,-18.8750){\line(1,-1){1.2500}}
\put(8.1250,-18.8750){\line(-1,-1){1.2500}}
\put(8.1250,-18.8750){\line(-1,1){1.2500}}
\put(5.6250,-20.1250){\makebox(2.5000,2.5000)[c]{\shortstack[c]{
\small Check\\
\small (B,C) }}} \put(5.6250,-18.5000){\makebox(0,0)[rt]{Violated}}
\put(8.1250,-18.5000){\makebox(0,0)[lt]{OK}}
\put(5.6250,-18.8750){\line(-1,0){3.0000}}
\put(2.6250,-18.8750){\line(0,1){5.0000}}
\put(8.1250,-18.8750){\line(1,0){2.0000}}
\put(10.1250,-18.8750){\line(0,1){6.0000}}
\put(10.1250,-12.8750){\line(1,0){1.0000}}
\put(11.1250,-12.8750){\vector(1,0){1.0000}}
\put(12.1250,-12.8750){\line(1,1){1.2500}}
\put(12.1250,-12.8750){\line(1,-1){1.2500}}
\put(14.6250,-12.8750){\line(-1,-1){1.2500}}
\put(14.6250,-12.8750){\line(-1,1){1.2500}}
\put(12.1250,-14.1250){\makebox(2.5000,2.5000)[c]{
 \shortstack[c]{
 \small $\|d\|<$\\ \small $10^{-\tol}$?}
 }}
\put(13.7500,-11.6250){\makebox(0,0)[lt]{No}}
\put(13.7500,-14.1250){\makebox(0,0)[lb]{Yes}}
\put(13.3750,-14.1250){\vector(0,-1){1.0000}}
\put(12.3750,-16.6250){\framebox(2.0000,1.5000)[c]{\shortstack[c]{
\small Predict\\
\small to t=0 }}} \put(13.3750,-16.6250){\vector(0,-1){1.0000}}
\put(11.3750,-19.6250){\line(1,1){2.0000}}
\put(11.3750,-19.6250){\line(1,-1){2.0000}}
\put(15.3750,-19.6250){\line(-1,-1){2.0000}}
\put(15.3750,-19.6250){\line(-1,1){2.0000}}
\put(11.3750,-21.6250){\makebox(4.0000,4.0000)[c]{\shortstack[c]{
\small Successive\\
\small predictions agree\\
\small to $10^{-\tol}$? }}}
\put(11.3750,-19.0250){\makebox(0,0)[rt]{No}}
\put(13.9750,-21.6250){\makebox(0,0)[lb]{Yes}}
\put(13.3750,-21.6250){\vector(0,-1){1.0000}}
\put(12.3750,-24.1250){\framebox(2.0000,1.5000)[c]{\shortstack[c]{
\small Return\\
\small endpoint }}} \put(11.3750,-19.6250){\line(-1,0){1.0000}}
\put(10.3750,-19.6250){\line(0,-1){2.0000}}
\put(10.3750,-21.6250){\vector(-1,0){1.0000}}
\put(6.3750,-22.6250){\framebox(3.0000,2.0000)[c]{\shortstack[c]{
\small Consider\\
\small increasing\\
\small Steplength }}} \put(6.3750,-21.6250){\vector(-1,0){1.0000}}
\put(2.3750,-22.6250){\framebox(3.0000,2.0000)[c]{\shortstack[c]{
\small Consider\\
\small decreasing\\
\small precision }}} \put(2.3750,-21.6250){\line(-1,0){1.0000}}
\put(1.3750,-21.6250){\line(0,1){14.5000}}
\put(1.3750,-7.1250){\line(1,0){2.0000}}
\put(13.3750,-11.6250){\vector(0,1){1.0000}}
\put(11.8750,-9.1250){\line(1,1){1.5000}}
\put(11.8750,-9.1250){\line(1,-1){1.5000}}
\put(14.8750,-9.1250){\line(-1,-1){1.5000}}
\put(14.8750,-9.1250){\line(-1,1){1.5000}}
\put(11.8750,-10.6250){\makebox(3.0000,3.0000)[c]{\shortstack[c]{
\small Iterations\\
\small remaining? }}} \put(13.8250,-7.6250){\makebox(0,0)[lt]{No}}
\put(11.8750,-8.6750){\makebox(0,0)[rt]{Yes}}
\put(11.8750,-9.1250){\line(-1,0){2.2500}}
\put(9.6250,-9.1250){\line(0,-1){6.2500}}
\put(9.6250,-15.3750){\vector(-1,0){1.5000}}
\put(13.3750,-7.6250){\line(0,1){3.7500}}
\put(13.3750,-3.8750){\vector(0,1){1.0000}}
\put(12.1250,-2.8750){\framebox(2.5000,1.5000)[c]{\shortstack[c]{
\small Decrease\\
\small Steplength }}} \put(12.1250,-2.1250){\vector(-1,0){5.2500}}
\end{picture}
\caption{Stepwise adaptive precision path tracking}
\end{figure}

The algorithm does not attempt corrections at $t=0$.  This is
because in our applications the target system $F(z)=H(z,0)$ often
has singular solutions.  It is safer to sample the incoming path
while it is still nonsingular and predict to $t=0$ based on these
samples.  In this situation, it helps to employ a more sophisticated
predictor than Euler's method.  For example, endgames that estimate
the winding number of the root and use it to compute a fractional
power series can be very effective~\cite{MSW91,MSW92b}.

\subsection{Error estimates}\label{Sec:ErrorEstimates}

To use the foregoing procedures, we need the function evaluation
error, $\psi$, and the errors contributing to $E$, namely, $\sE$ and
$\phi$.  There is a trade-off between using rigorously safe bounds
for highest reliability or using less stringent figures reflecting
typical behavior to avoid the overuse of high precision.  Rough
figures are acceptable as this is just a means of setting the
precision.  Also, a user of the path tracker will not usually wish
to expend a lot of effort in developing error bounds.

A rigorous and automated way of establishing error bounds is to use
interval arithmetic.  Following that approach, one may wish go all
the way and use interval techniques to obtain a path tracker with
fully rigorous step length control, as in \cite{KX94}.  However,
this can be expensive, due partially to the cost of interval
arithmetic but more significantly due to the cost of
overconservative error bounds, which slow the algorithm's progress
by driving the step size smaller than necessary.  Still, when
rigorous results are desired, it may be worth the cost.  The method
of \cite{KX94} does not explicitly include adaptive precision, so
something along the lines discussed here could be useful in
modifying that approach.

Instead of using interval methods, we may approximate errors by
accumulating their effects across successive operations.  Suppose
the program to evaluate $F(z)$ has been parsed into a straight line
program, that is, a sequence of unary and binary operations free of
branches or loops. Suppose that at some intermediate stage of
computation, we have computed a value $\hat a$ for a real number
$a$, such that $\hat a$ lies between $(1-u)a$ and $(1+u)a$. (For a
floating point complex number, this applies to both the real and
imaginary parts.) Let's use the shorthand $\hat a = (1\pm u)a$ to
mean this entire interval. If $\hat a = (1\pm u_a)a$ and $\hat
b=(1\pm u_b)b$, the product $c=ab$ is computed in floating point
with unit round-off $\ubar$ as $\hat c \approx
ab(1\pm\max[\ubar,(u_a+u_b)])$, where the quadratic round-off term
$u_au_b$ is neglected.  The absolute error $\hat c - c$ thus has the
approximate bound $\pm\max[\ubar,(u_a+u_b)]|a| |b|$. Similarly,
$a+b$ is computed as $a+b\pm u_a a \pm u_b b$ which has an absolute
error bounded by $\pm\max[\ubar|a+b|, (u_a|a|+u_b|b|)]$. Using just
these relations, an error bound for any straight-line polynomial
function can be calculated in parallel with the function evaluation
itself. Similar relations can be developed for any smooth elementary
function, such as the basic trigonometric functions. Assuming that
the inputs to the function, including both the input variables $z$
and any internal parameters of the function, are all known either
exactly or with relative round-off error $\ubar$, the output of the
error analysis is $\Phi(z)$ such that the error in the computed
value $\hat F(z)$ is $\|\hat F(z)-F(z)\|=\psi(z)=\Psi(z)\ubar$. When
the result is rounded off to a possibly lower precision $u$, the
total error becomes the form shown in Eq.~\ref{Eq:errorF}.

It is important to note that the error in the function depends on
the error in its parameters.  For example, consider the simple
function $f(x)=x^2-1/3$.  If this is rounded off to $g(x)=x^2-0.333$
before we use high precision to solve $g(x)=0$, we will obtain an
accurate value of $\sqrt{0.333}$ but we will never get an accurate
value of $1/\sqrt{3}$.  Although this is an obvious observation, it
can easily be forgotten in passing a homotopy function from some
application to the adaptive precision path tracking algorithm.  If
coefficients in the function are frozen at fixed precision, the
algorithm tracks the solutions of the frozen function, not the exact
problem that was intended.  Whether the difference is significant
depends on the nature of the application and the sensitivity of the
function.

While $\psi$ and $\phi$ concern the errors in evaluating the
function and its Jacobian, the factor $\sE$ concerns the stability
of the linear solve.   Round-off errors can accumulate through each
stage of elimination.  When Gaussian elimination with partial
pivoting is used, the worst-case error bound grows as $\sE=2^n$ for
solving an $n\times n$ system \cite{DemmelBook}.
However, as indicated in \cite{DemmelBook}, $\sE$ rarely exceeds $n$ with
the average case around $n^{\frac{2}{3}}$ or $n^{\frac{1}{2}}$.  Setting
$\sE = n^2$ should therefore be sufficient for almost all cases.

\subsection{Application to polynomial systems}\label{Sec:PolySystems}

To avoid program complexity and save computation time, it is
preferable not to perform a full error analysis of the type just
described.  In many cases a rough analysis is sufficient and easily
derived.  This is indeed possible for the case of most interest to
us: polynomial systems.

Suppose $h:\bC^{n+1}\rightarrow\bC$ is a degree $D$ homogeneous
polynomial
\begin{equation}\label{Eq:HomogH}
  h(x)=h(x_0,x_1,\ldots,x_n)=\sum_{i\in\sI} c_i x_0^{d_{0i}}\cdots
  x_n^{d_{ni}}, \qquad \sum_{j=0}^n d_{ji}=D, \hbox{ for all
  $i\in\sI$},
\end{equation}
where $\sI$ is just an index set for the coefficients $c_i$. Since
$h$ is homogeneous, $h(\lambda x)=\lambda^Dh(x)$, so if $h(x)=0$,
then also $h(\lambda x)=0$.  Consequently, the solution set of a
system of homogeneous polynomials can be said to lie in projective
space $\bP^n$, the set of lines through the origin in $\bC^{n+1}$.
Similarly, the solutions of multihomogeneous polynomials lie in a
cross product of projective spaces, see \cite{SoWa}.

Any inhomogeneous polynomial $g(z_1,\ldots,z_n)$ can be easily
homogenized to obtain a related function $G(x_0,x_1,\ldots,x_n)$,
with
\[
 G(1,x_1,\ldots,x_n)=g(x_1,\ldots,x_n).
\]
Hence, for any solution $z_*$ of $g(z)=0$ there is a corresponding
solution $x_*=(1,z_*)$ of $G(x)=0$.  One advantage of homogenization
is that we can re-scale any solution $x_*$ of $G(x)=0$ to make
$\|x_*\|=1$, which often helps numerical conditioning.

Error bounds for homogeneous polynomials can be estimated easily. If
we rescale $x$ so that the maximum entry in $(x_0,\ldots,x_n)$ has
magnitude 1, then the error in evaluating the degree $D$ homogeneous
polynomial $h(x)$ as in Eq.~\ref{Eq:HomogH} is approximately
\begin{equation}\label{Eq:HomogPsi}
    \psi(x,u) \approx uD\sum_{i\in\sI} |c_i|,
       \hbox{ i.e., }\Psi=D\sum_{i\in\sI} |c_i|.
\end{equation}
Similarly, the derivatives have an approximate error bound of
\begin{equation}\label{Eq:HomogPhi}
    \phi(x,u) \approx uD(D-1)\sum_{i\in\sI} |c_i|,
       \hbox{ i.e., }\Phi=D(D-1)\sum_{i\in\sI} |c_i|.
\end{equation}

At first glance, it may seem that errors can be reduced by simply
scaling the functions, and thereby scaling their coefficients, by
some small factor. But $\|J_*^{-1}\|$ will scale oppositely, so the
error predicted by Eq.~\ref{Eq:FinalError} is unchanged.

\section{Computational results}\label{Sec:Results}

This section contains a brief discussion of the implementation
details for multiprecision arithmetic and for evaluating the rules
for adapting precision.  Then we discuss the results of applying the
adaptive precision path tracker to three example polynomial systems.

\subsection{Implementation details}\label{Sec:Details}

Bertini is a software package for computation in numerical algebraic
geometry currently under development by the authors and with some
early work by Christopher Monico.  Bertini is written in the C
programming language and makes use of straight-line programs for the
representation, evaluation, and differentiation of polynomials.  All
the examples discussed here were run using an unreleased version of
Bertini on an Opteron 250 processor running Linux.

The adaptation rules, A, B, and C (or \Cp), leave some choices open
to the final implementation.  For the runs reported here, we chose
to evaluate function residuals to the same precision as the
computation of Newton corrections, so rule C applied, not rule \Cp.
Also, in rules A and B, we chose to use $\sE=n^2$, where $n$ is the
number of variables, which is somewhat conservative for typical
cases but underestimates the worst pathological cases. (See
Section~\ref{Sec:ErrorEstimates} for more on this issue.) The rules
require formulas for evaluating the error bounds $\phi(x,u)=\Phi u$
and $\psi(x,u)=\Psi u$.  These are problem dependent, so we report
our choices for each of the example problems below.

To adaptively change precision, Bertini relies on the open source
MPFR library for multiprecision support.  Bertini has data types and
functions for regular precision (based on the IEEE ``double''
standard) and higher precision (using MPFR).  Although the program
would be simpler if MPFR data types and functions were used
exclusively, the standard double precision types and functions in C
are more efficient, so Bertini uses these whenever the adaptation
rules indicate that double precision is sufficient. Additional
details regarding the use of multiple precision may be found using
links from the Bertini website.  Since the use of adaptive precision
variables is highly implementation-specific, no other details are
described here.

MPFR requires the addition of precision to the mantissa in packets
of 32 bits. Since the discussion of the examples below involves both
binary and decimal digits, Table~\ref{Tab:Bits2Digits} shows how to
convert between the two.

\begin{table}
  \centering
  \begin{tabular}{|c|c|c|c|c|c|c|}
  \hline
  &IEEE single&IEEE double&\multicolumn{4}{c|}{MPFR}\\
  \hline
  bits&23&52&64&96&128&256\\
  \hline
  decimal digits&6&15&19&28&38&77\\
  \hline
  \end{tabular}
  \caption{Number of digits for
  mantissas at various levels of precision}\label{Tab:Bits2Digits}
\end{table}

\subsection{Behavior of adaptive precision near a
singularity}\label{Sec:Griewank}

Consider the polynomial system of Griewank and Osborne \cite{GO83},
$$
f = \left[\begin{array}{l}
\frac{29}{16}z_1^3 - 2z_1z_2\\
z_2 - z_1^2\\
\end{array}
\right].
$$
This system is interesting because Newton's method diverges for any
initial guess near the triple root at the origin. A homotopy of the
form
\[
  h(z,t) = tg(z) + (1-t)f(z),
\]
where $g(z)$ is the following start system, constructed as
described in \cite{SoWa}:
$$
g1=\begin{array}{l}
((-0.74924187 + 0.13780686i)z_2 + (+0.18480353 - 0.41277609i)H_2)*\\
((-0.75689854 - 0.14979830i)z_1 + (-0.85948442 + 0.60841378i)H_1)*\\
((+0.63572306 - 0.62817501i)z_1 + (-0.23366512 - 0.46870314i)H_1)*\\
((+0.86102153 + 0.27872286i)z_1 + (-0.29470257 + 0.33646578i)H_1),\\
\end{array}
$$
$$
g2=\begin{array}{l}
((+0.35642681 + 0.94511728i)z_2 + (0.61051543 + 0.76031375i)H_2)*\\
((-0.84353895 + 0.93981958i)z_1 + (0.57266034 + 0.80575085i)H_1)*\\
((-0.13349728 - 0.51170231i)z_1 + (0.42999170 + 0.98290700i)H_1),\\
\end{array}
$$
with $H_1$ and $H_2$ the extra variables added to the variable
groups $\{z_1\}$ and $\{z_2\}$, respectively.
In the process of homogenizing, two linear polynomials (one per variable
group) are added to the system and do not depend on $t$.  In this
case, those polynomials are:
\[\begin{array}{l}
(-0.42423834 + 0.84693089i)z_1 + H_1 - 0.71988539 + 0.59651665i,\\
(+0.30408917 + 0.78336869i)z_2 + H_2 + 0.35005211 - 0.52159537i.\\
\end{array}
\]
This homotopy has three
paths converging to the origin as $t\rightarrow 0$.  These paths
remain nonsingular for $t\in(0,1]$, so it is interesting to see how
the adaptive precision algorithm behaves as $t$ approaches zero.

Using the prescription of Eqs.~\ref{Eq:HomogPsi},\ref{Eq:HomogPhi}
to bound errors, we take $D=3$ and $\sum_{i\in\sI} |c_i|\approx 4$,
so $\Psi=12$ and $\Phi=24$. (This is actually quite conservative,
because a more realistic bound would be $\Psi=12||z||$ and $||z||$
is approaching zero.)  We set the safety digits
$\sigma_1=\sigma_2=0$, as these serve only to force precision to
increase at a slightly larger value of $t$.

We set the desired accuracy to $\tau=8$ for $t\in[0.1,1]$ and
increased it to $\tau=12$ digits thereafter. To watch the behavior
of the algorithm for very small values of $t$, we turned off the
usual stopping criterion and instead simply ran the path to
$t=10^{-30}$. That is, in the flowchart of
Figure 2, after a successful correction
step, the algorithm was modified to always loop back to step ahead
in $t$ until $t=10^{-30}$. Since for small $t$ and for $||g||\approx
1$, the homotopy path has $||z||\approx |t|^{1/3}$, all three roots
heading to the origin were within $10^{-10}$ of the origin at the
end of tracking. With an accurate predictor, such as a fractional
power series endgame \cite{MSW92b}, the solution at $t=0$ can be
computed to an accuracy of $10^{-12}$ using only double precision,
but the purpose of this example is to show that the adaptive
precision tracking algorithm succeeds beyond the point where double
precision would fail.

The result is shown in Figure 3. We plot the
right-hand side of rule C, as this rule determines the increases in
precision for this problem.  The jump in this value at $t=0.1$ is
due to our prescribed increase in $\tau$ at that point.  Since the
path is heading to the origin, if we used the less conservative
error estimate of $\psi(z,u)=12||z||u$, increases in precision would
be delayed somewhat, as rule C would be shifted down by
approximately $\log_{10}||z||\approx (1/3)\log_{10} t$.

\begin{figure}\label{Fig:GriewankPlot}
\centerline{
\epsfxsize=3in
\epsfysize=3in
\epsfbox{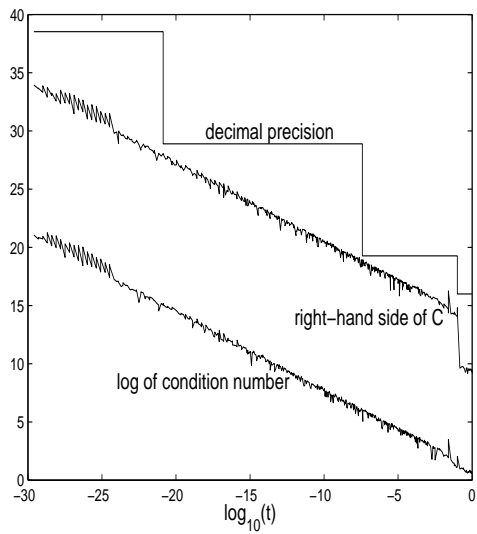}
}
\caption{\small Right-hand side of C, condition number, and decimal precision versus $\log(t)$}
\end{figure}

In cases where multiple roots are clustered close together or even
approaching the same endpoint, as in this example, the path tracking
tolerance must be kept substantially smaller than the distance
between roots to avoid jumping between paths.  Rather than
prescribing $\tau$ versus $t$, as we did in this example, it would
be preferable to automatically adjust $\tau$ as needed.  We postpone
this for future research, but note that the rules given here for
adjusting precision should still be effective when $\tau$ is
adjusted automatically.

\subsection{Behavior of adaptive precision under tighter tolerances}\label{Sec:tol_example}

To illustrate the effect of tightening the  tracking tolerance
(i.e., increasing $\tau$) on adaptive precision path tracking, we
will consider a polynomial system coming from chemistry.  To
determine chemical equilibria, one may pass from a set of reaction
and conservation equations to polynomials, as described in
\cite{Mor87,SoWa}.  One such polynomial system, discussed in
\cite{MSW92a}, is the following:
$$
f=\left[\begin{array}{l}
14z_1^2 + 6z_1z_2 + 5z_1 - 72z_2^2 - 18z_2 - 850z_3
+ 0.000000002\\
0.5z_1z_2^2 + 0.01z_1z_2 + 0.13z_2^2 + 0.04z_2 -
40000\\
0.03z_1z_3 + 0.04z_3 - 850\\
\end{array}\right].
$$
As in the previous example, this system was homogenized, although
this time with only one variable group, so there is just a single
homogenizing coordinate. Using a compatible total-degree linear
product start system, that is, one whose polynomials have degrees 2,
3, and 2, respectively, results in a homotopy with 12 paths. The
solution set of this system consists of eight regular finite solutions, two
of which have some coordinates of size around $10^5$.
These eight finite solutions are provided in Table~\ref{Tab:ChemSysSolns}.

\begin{table}
  \centering
  \tiny
  \renewcommand{\arraystretch}{1.25}
  \begin{tabular}{|c|r@{\hspace{4pt}}c@{\hspace{4pt}}l|r@{\hspace{4pt}}c@{\hspace{4pt}}l|}
  \hline
  $H_1$&2.15811678208e$-$03&$ - $& 2.32076062821e$-$03*i&7.75265879929e$-$03&$ + $& 5.61530748382e$-$03*i\\
  $z_1$&1.21933862567e$-$01&$ + $& 4.02115643024e$-$01*i&1.26177608967e$-$01&$ - $& 1.26295173168e+00*i\\
  $z_2$&-2.29688938707e$-$02&$ - $& 7.44021609426e$-$02*i&-2.45549175888e$-$02&$ + $& 2.33918398619e$-$01*i\\
  $z_3$&-6.62622511387e$-$01&$ - $& 1.55538216233e$-$00*i&-1.87267506123e+00&$ + $& 8.48441541195e$-$01*i\\
  \hline
  $H_1$&-2.54171295092e$-$03&$ + $& 1.43777404446e$-$03*i&6.32152240723e$-$03&$ + $& 1.64022773970e$-$03*i\\
  $z_1$&-3.34864497185e$-$01&$ + $& 1.89423218369e$-$01*i&-2.20546409488e$-$01&$ - $& 7.88700342178e$-$01*i\\
  $z_2$&6.25969991088e$-$02&$ - $& 3.54093238711e$-$02*i&-4.39498685300e$-$02&$ - $& 1.59117743373e$-$01*i\\
  $z_3$&-5.41138529778e$-$01&$ + $& 3.06106507778e$-$01*i&-1.03394901752e+00&$ + $& 1.06381058693e+00*i\\
  \hline
  $H_1$&-2.64917471213e$-$05&$ + $& 5.83090377404e$-$06*i&2.72428081371e$-$03&$ - $& 2.22348561510e$-$03*i\\
  $z_1$&1.23749450722e$-$05&$ - $& 2.72846568805e$-$06*i&6.93769305944e$-$02&$ + $& 4.35467392206e$-$01*i\\
  $z_2$&-3.62126436085e$-$03&$ - $& 1.64631735533e$-$02*i&1.42630599439e$-$02&$ + $& 8.77317115664e$-$02*i\\
  $z_3$&-8.66534113884e$-$01&$ + $& 1.90904229879e$-$01*i&-7.45011925697e$-$01&$ - $& 2.88085192442e$-$01*i\\
  \hline
  $H_1$&-2.37392215058e$-$03&$ + $& 9.07039153390e$-$04*i&-2.73688783636e$-$05&$ + $& 4.95136100653e$-$06*i\\
  $z_1$&-2.96180844307e$-$01&$ + $& 1.13166145980e$-$01*i&1.27865341710e$-$05&$ - $& 2.30844830185e$-$06*i\\
  $z_2$&-6.00257143378e$-$02&$ + $& 2.29349024594e$-$02*i&3.07919899933e$-$03&$ + $& 1.70073434711e$-$02*i\\
  $z_3$&-5.33404946327e$-$01&$ + $& 2.03805834055e$-$01*i&-8.95294964314e$-$01&$ + $& 1.61788761616e$-$01*i\\
  \hline
  \end{tabular}
  \normalsize
  \caption{The solutions of the chemical system}\label{Tab:ChemSysSolns}
\end{table}

Due to the poor scaling of the problem, the error bounds were set to
$\Psi=120,000$ and $\Phi=240,000$.  As opposed to the previous example, the usual
stopping criterion described in Section \ref{Sec:Newton} was
employed.  No endgame was used, though, since the use of endgames
speeds convergence.

Tracking to a final tolerance of $10^{-8}$ using fixed
regular (IEEE double) precision, all eight finite solutions were
discovered, including the two having large coordinates.
In tightening the tolerance to
$10^{-12}$, however, the linear algebra involved with
path tracking broke down for the paths leading to these
two large solutions.  This resulted in step failures and,
ultimately, path failure for those two paths.  However, by
increasing precision to 96 bits, all eight
regular solutions were again found.

Using adaptive precision with the number of safety digits set to 0,
1, or 2, the six finite solutions of moderate size required only one
precision increase (to 64 bits from 52 bits).  This increase
occurred at the very end of the path.  The two large finite
solutions each needed 96 bits of precision (as expected).

Since the initial run in double precision succeeded on the six
moderate paths, the adaptive method's increase to 64 bits at the end
of these paths is not necessary.  For safety, the adaptive rules are
designed to be conservative, so some extra computational cost is to
be expected.

\subsection{A class of univariate polynomials}\label{Sec:Cheby}

The Chebyshev polynomials have been studied extensively and are
known to have many interesting properties \cite{DemmelBook}. There
is one Chebyshev polynomial in each degree, and scaled such that the
leading coefficient is 1, they may be defined recursively by
$$
\begin{array}{l}
T_0(x) := 2,\\
T_1(x) := x,\text{ and}\\
T_i(x) := xT_{i-1}(x) - T_{i-2}(x)/4, \text{ for } i \ge 2.
\end{array}
$$
The solutions of $T_{n+1}(x)$ are then given by
$$
\cos\left(\frac{\left(2n+1-2k\right)\pi}{2n+2}\right),
$$
for $k=0,1,...,n$.

Several of these polynomials, with degrees ranging from 10 to 300,
were solved in Bertini using the step-adaptive precision path
tracking method developed above.  The systems were run
with a linear homotopy with a total degree start system and without
homogenizing.  In each case, bounds on $\Phi$ and $\Psi$ were set
almost exactly.  In particular, using the estimates of $\Phi$ and
$\Psi$ described in Section~\ref{Sec:PolySystems}, $D$ was chosen
exactly, and $\sum_{i\in\sI} |c_i|$ was set to the values given in
Table~\ref{Tab:AbsValsOfCoeffs}.  In each case, both safety
digit settings were set to 4 and $\tau$ was set to 10.

\begin{table}
  \centering
  \begin{tabular}{|c|c|}
  \hline
  Degree&Estimate\\
  \hline
  10&7\\
  15&17\\
  20&44\\
  25&111\\
  50&12300\\
  100&$1.5\cdot10^8$\\
  150&$1.9\cdot10^{12}$\\
  200&$2.3\cdot10^{16}$\\
  250&$2.8\cdot10^{20}$\\
  300&$3.4\cdot10^{24}$\\
  \hline
  \end{tabular}
  \caption{Estimates of $\sum_{i\in\sI} |c_i|$ for various degrees}\label{Tab:AbsValsOfCoeffs}
\end{table}

The level of precision used by the step-adaptive precision path
tracking algorithm developed above was degree- and path-dependent
for the Chebyshev polynomials, although all paths for a given degree
needed approximately the same level of precision.  In each degree
that was considered, the path ending nearest $1.0$ was one of the
paths needing the highest level of precision for that degree. (This
occurs because the spacing between the roots is smallest near
$\pm1$.) The levels of precision used for that path are displayed in
Figure 4. It should be noted that for every degree considered, a
complete solution set with all solutions correct to at least 10
digits was found.

\begin{figure}[h]\label{Fig:ChebyPrec}
\centerline{
\epsfxsize=3in
\epsfysize=3in
\epsfbox{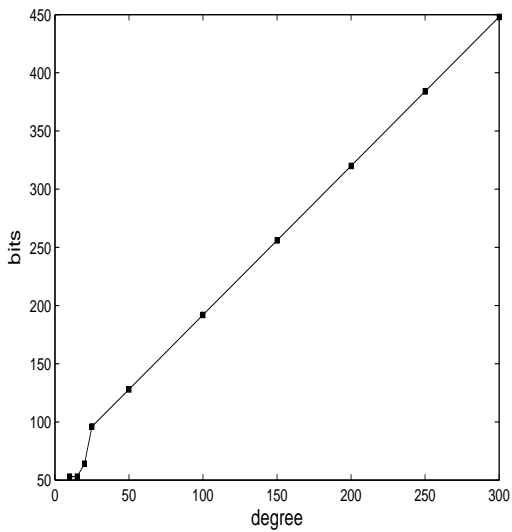}
}
\caption{\small Number of bits needed for Chebyshev polynomials of various degrees}
\end{figure}

We note that to solve the high degree Chebyshev polynomials, a small
initial step size was required to get the path tracker started. With
too large an initial step, the predicted point was so far from the
path that the adaptive precision rules increased precision to an
unreasonable level without ever exiting the corrector loop.  As
diagrammed in Figure~2, the algorithm must exit the corrector loop
before a decrease in step length can be triggered.  Various ad hoc
schemes could detect and recover from this type of error, but we
would prefer a step size control method based on an analysis of the
predictor. For the moment, we defer this for future work.

\section{Discussion}

On some problems, endgames can speed convergence to the point that
singular endpoints can be estimated accurately in double precision.
It should be noted that is not generally the case. Without enough
precision, the ``endgame operating zone'' is empty
\cite{MSW92b,SoWa}. Likewise, endgames based on deflating the system
to derive a related nonsingular one \cite{LVZ} may need higher than
double precision to make a correct decision on the rank of the
Jacobian at each stage of deflation \cite{LZ}. Moreover, if some
other sort of singularity is encountered during path tracking, away
from $t=0$, endgames will not be useful while adaptive precision
will be. In the case of tight final tolerances or endpoints of paths
having high multiplicity, endgames will again need assistance from
higher (and therefore adaptive) precision methods.  Conversely, high
precision is expensive and floating point precision can never be
made truly infinite, so to get the most out of whatever precision
one uses, endgames are indispensable.

The theory going into this new adaptive precision method revolves
around Newton's method or corrector methods in general.  However,
corrector methods are only one half of basic path tracking.  A
careful study of predictor methods is certainly warranted.  The use
of different predictor schemes, e.g., Adams-Bashforth rather than
Euler, is well worth considering.  A careful analysis of the
predictor might be combined with the convergence criteria of the
corrector to automatically determine a safe step length in place of
the trial-and-error step length adaptation method we have used here.
This might give an efficient alternative to \cite{KX94}, which
presents a rigorous step length control algorithm based on interval
arithmetic.

Another open question, discussed briefly in
Section~\ref{Sec:Griewank}, is the issue of adaptively changing the
tracking tolerance in response to close approaches between paths.

\end{document}